\documentclass[a4paper,12pt]{article}
\usepackage{amsmath}
\usepackage{amsfonts}
\usepackage{amssymb}

\newtheorem{remark}{Remark}[section]
\newtheorem{theorem}{Theorem} [section]
\newtheorem{lemma}{Lemma}[section]

\newtheorem{proposition}{Proposition}[section]
\newcommand\qed{\hfill\kern1pt \vbox{\hrule height 0.6pt \hbox{\vrule width 0.6pt \vbox{\vskip 6pt}  \hskip 3pt \vrule width 1.3pt} \hrule depth 1.3pt}     \kern1pt\medskip\noindent}

\newcommand{\be}{\begin{equation}} 
\newcommand{\ee}{\end{equation}}
\newcommand{\bea}{\begin{eqnarray}} 
\newcommand{\eea}{\end{eqnarray}}
\newcommand{\bean}{\begin{eqnarray*}} 
\newcommand{\eean}{\end{eqnarray*}}
\newcommand\R{{\mathbb R}}

\renewcommand\S{{\mathcal S}} 
\newcommand\X{{\mathcal X}}  
\def\f{\varphi}
\def\d#1{\, {\rm d}#1}

\def\Proof{{\smallskip\noindent{\em Proof. }}}

\begin{document}
     \baselineskip=17pt
\title{
On the parabolic-elliptic limit \\of the
doubly parabolic Keller--Segel system\\
 modelling chemotaxis}  
\author{Piotr Biler$^1$ and Lorenzo Brandolese$^2$\\ 
\small{ $^1$ Instytut Matematyczny, Uniwersytet Wroc{\l}awski,}\\ 
\small{ pl. Grunwaldzki 2/4, 50--384 Wroc{\l}aw, POLAND}\\ 
\small{\tt Piotr.Biler@math.uni.wroc.pl}\\ 
\small{ $^2$  Universit\'e de Lyon, Universit\'e Lyon 1,}\\
\small{Institut Camille Jordan, CNRS UMR 5208,}\\
\small{43 bd. du 11 Novembre,  69622 Villeurbanne Cedex, FRANCE }\\ 
\small {\tt brandolese@math.univ-lyon1.fr}}

\date{\today}

\maketitle

\begin{center}
\emph{Dedicated to the memory of  Andrzej Hulanicki}
\end{center}
\vskip0.5cm

\begin{abstract}

We establish new convergence results, in strong  topologies,
for solutions
of the parabolic-parabolic Keller--Segel system in the plane, to the corresponding solutions of the
parabolic-elliptic model, as a~physical parameter goes to zero.
Our main tools are suitable space-time estimates,
implying the global existence of slowly decaying (in general, nonintegrable)
solutions for these models, under a natural smallness assumption.
\end{abstract}

\section{Introduction}

We consider two related nonlinear parabolic systems which are frequently used as models for a description of chemotactic phenomena, including the aggregation of microorganisms caused by a chemoattractant, i.e. a  chemical whose concentration gradient governs the oriented movement of those microorganisms. The parabolic character of the systems comes from the diffusion described by the Laplacians. A version of the system (PE) below is also used in astrophysics as a model of the evolution of a cloud of self-gravitating particles in the mean field approximation.

The first one is the classical parabolic-elliptic Keller--Segel system
\begin{equation*}
\left\{
\begin{aligned}
&u_t =\Delta u-\nabla\cdot(u\nabla \varphi),\\
&\Delta \varphi+u=0,\\
&u(0)=u_0.\;
\end{aligned}
\right.
\qquad x\in\R^2,\ t>0,
\tag{PE}
\end{equation*} 
Here, $u=u(x,t)$, $\f=\f(x,t)$ are either functions or suitable (tempered) distributions. 
When $u\ge 0$, $\f\ge 0$, they may be interpreted as concentrations (densities) of microorganisms and chemicals, respectively. 

The second one is the parabolic-parabolic system 
\begin{equation}
\left\{
\begin{aligned}
&u_t =\Delta u-\nabla\cdot(u\nabla \varphi),\\
&\tau \varphi_t=\Delta \varphi+u,\\
&u(0)=u_0,\ \  \varphi(0)=0,
\end{aligned}
\right.
\qquad x\in\R^2,\ t>0,
\tag{PP}
\end{equation}
where $\tau>0$ is a fixed parameter. Each of the models can be considered as a~single nonlinear parabolic equation for $u$ with a \emph{nonlocal} (either in $x$ or in $(x,t)$) nonlinearity since the term $\nabla \f$ can be expressed as a linear integral operator acting on $u$. In the latter model, the variations of the concentration  $\f$ are governed by the linear nonhomogeneous heat equation, and therefore are slower than in the former system, where the response of $\f$ to the variations of $u$ is instantaneous, and described by the integral operator $(-\Delta)^{-1}$ whose kernel 
has a~singularity. Thus, one may expect that the evolution described by (PE) might be faster than that for (PP), especially for large values of $\tau$ when the diffusion of $\f$ is rather slow compared to that of $u$. Moreover, the nonlinear effects for (PE) should manifest themselves faster than for (PP). 
\medskip

The theory of the system (PE) is relatively well developed, in particular when this is studied in a bounded domain in $\R^d$, $d=1,\, 2,\, 3$, with the homogeneous Neumann conditions for $u$ and $\f$ at the boundary of the domain. One of the most intriguing properties of (PE) considered for \emph{positive} and \emph{integrable} solutions $u$ in $d=2$ case is the existence of a threshold value $8\pi$ of mass $M\equiv \int u(x,t)\d{x}$, see the pioneering work \cite{JL} and \cite{B3,B-AMSA}. Namely, if $u_0\ge 0$ is such that $\int u_0(x)\d{x}>8\pi$, then any regular, positive  solution $u$ of (PE) cannot be global in time. We refer the reader for a fine description of the asymptotic behaviour of integrable solutions of (PE) in the subcritical case $M<8\pi$ to \cite{BDP} and to~\cite{BCM} for the limit case~$M=8\pi$.
See also \cite{BKLN1,BKLN2} in the radially symmetric case.
 The higher dimensional versions of (PE) have been also extensively studied, cf., e.g., \cite{B-SM,BCGK,BB}, and \cite{B3} for blow up phenomena. 
\medskip

The doubly parabolic system (PP) has been a bit less studied. For instance, it is known that if for the initial data $u_0$ one has $M<8\pi$,  then positive solutions are global in time, see \cite{CC} in the case of a bounded planar domain, and also \cite[Theorem 5]{B-AMSA}. However, it is not known whether $M\le 8\pi$ is, in general, a \emph{necessary} condition for the existence of global in time solutions. That is,  it is not known whether the blow up occurs for solutions, except for a specific example in \cite{HV} of a~particular blowing up solution for a system close to (PP). Even, it is an open question what is the exact range of $M$ guaranteeing the existence of integrable \emph{self-similar solutions}.  For the system (PE) it is proved that  $M\in[0,8\pi)$, and the self-similar solutions (unique for a given $M\in[0,8\pi)$) describe the generic asymptotic behaviour of global in time, positive  and integrable solutions of (PE). Concerning (PP), it is known that $M<M(\tau)$ with $M(\tau)$ linear in $\tau$, is a necessary condition for the existence of self-similar solutions, cf. \cite{B-BCP}. For (PP) with small $M$ such special solutions are also important in the study of space-time decay of general solutions, see \cite{N}. The analysis if any $M>8\pi$ may correspond to a self-similar solution is under way, see \cite{BCD}.
For a different point of view about self-similar solutions for higher dimensional models of (PP), see also~\cite{KozS08}. 
Usual proofs of a blow up for (PE) involve calculations of moments of a solution and then symmetrization, cf. \cite{BDP,B3}. These methods seem do not work for (PP), hence another approach is needed to show a~blow up for that system.
For a numerical insight on blow-up issues we refer, e.g.,  to~\cite{Filbet06}.
\medskip

A nice result in \cite{AR} shows that the solutions of the systems (PP) and (PE) enjoy a kind of stability property as $\tau\searrow 0$: solutions of (PP) converge in a suitable sense to those of (PE). It had been an old question raised by J. J. L. Vel\'azquez and D. Wrzosek, recently solved in \cite{AR}. However, this result obtained for suitably small solutions in quite a big functional space of pseudomeasures, gives no indication on the behaviour of  possible (``large'') blowing up solutions. 
\medskip

The solvability of the systems (PE), (PP) has been studied in various classes of functions and distributions, like Lebesgue, Morrey, Besov, etc., with an immediate motivation to include the \emph{a~priori} strongest possible critical singularities of either solutions or initial data which appear to be point measures in the two-dimensional case and the multiples of $|x|^{-2}$ function in the higher dimensional case. In particular, ``vast'' functional spaces suitable for analysis of the two-dimensional systems include measure and pseudomeasure spaces, cf. \cite{B-SM,B-AMSA,BCGK,AR}. 
\medskip

We show in this paper a result on the existence of (in general, nonintegrable) solutions in a
class~$\X$ of functions with natural space-time decay properties,
see Theorem~\ref{theorem1} and \ref{theorem2}.
Here, the space~$E$ of admissible initial conditions also contains Dirac measures.  The corresponding solutions may be positive and ``large'' in the sense of their nonintegrability. Nevertheless, they are defined globally in time.
Unlike the paper \cite{AR}, we work in $(x,t)$ space, when \cite{AR} has dealt with the Fourier variables $\xi$,
cf. the formulation \eqref{F} below. 
Such results are obtained by an extension and refinement of techniques used in \cite{BB} for (PE) in higher dimensions, but neither for (PP),
nor in the two-dimensional case of (PE) which often requires a specific treatement.
Moreover, the function spaces that we employ here 
allow us to deal with data that can be more singular than those considered in~\cite{BB}.
The spaces $\X$ and~$E$ defined in the next section are, in a sense, critical for that analysis, and have
been already considered, in slightly different forms,
e.g. in the studies of the Navier--Stokes system in \cite{CDW, M-2}. 

\medskip
Our main results are contained in Section~\ref{section4}, where we address the problem of the convergence as $\tau\searrow0$
of solutions $u^\tau$ of the system (PP) to the corresponding solutions $u$ of (PE),
in the space~$\X$, arising from small data in~$E$.
Mathematically, our stability result
is not included in, and does not imply, that of~\cite{AR}.
However, it seems to us that the use of the natural $(x,t)$ variables
provides a more immediate physical interpretation.
Furthermore, our method looks more flexible, and can be used
to prove the stability of the system with respect to stronger topologies.
For example, we establish also the convergence in the $L^\infty_t(L^1_x)$-norm
for data belonging to $E\cap L^1$, and in the $L^\infty_{t,x}\,$-norm for data in~$E\cap L^\infty$.
Motivated by~\cite{LR}, we will also address this issue in the
more general setting of~\emph{shift invariant spaces of local measures\/}.
The main difficulty for obtaining the convergence $u^\tau\to u$ in strong norms is that
$\nabla\varphi^\tau$ enjoys some kind of \emph{instability} as $\tau\to0$,
in particular in weighted spaces.

\medskip
Moreover, we give a nonexistence (blow up) result for solutions of (PP) in $\R^d$, $d\ge 1$, with the positive Fourier transform of $\widehat{u_0}$ in the spirit of \cite{MS}, see Theorem~\ref{theoremCC}. These are complex valued solutions with no straightforward physical/biological interpretation.
However, such a result tells us that there is no hope to prove the global existence of solutions to (PP) and similar models for arbitrarily large
data relying only on size estimates.

\bigskip
\section{The parabolic-elliptic system}

In order to study the systems (PE) and (PP) we introduce the Banach space
$\X$ of functions $u=u(x,t)$ and the Banach space $E$ 
of tempered distributions $u_0\in\mathcal{S}'(\R^2)$
by defining the norms
\begin{equation}
\label{e1}
\|u\|_\X=\hbox{ess\,sup}_{t>0,x\in\R^2} \, (t+|x|^2)\, |u(x,t)|,
\end{equation}
and
\begin{equation}
\|u_0\|_E=\|e^{t\Delta}u_0\|_\X.
\end{equation}
Here, $e^{t\Delta}$ denotes the heat semigroup defined by the Gaussian kernel $g_t$, $g_t(x)=\frac{1}{4\pi t}e^{-\frac{1}{4t}|x|^2}$. 
For example, the Dirac mass $u_0=\delta$ in $\R^2$ is an element of $E$.
Notice that, by the definition, $E$ is continuously embedded into the weak Hardy space $\mathcal{H}^1_w$,
which is the space consisting of all tempered distributions~$f$ such that
$\sup_{t>0} |e^{t\Delta}f|$ belongs to the Lorentz space $L^{1,\infty}(\R^2)$.
See \cite{M-1}.

\medskip

Let us define the bilinear form $B_0$ by 
\begin{equation} \label{B0}
B_0(u,v)(t)\equiv \int_0^t e^{(t-s)\Delta}\nabla \cdot(u\nabla(-\Delta)^{-1} v)(s)\d{s}.
\end{equation} 
Here, $(-\Delta)^{-1}$ is the convolution operator on functions defined on $\R^2$ with the kernel $K(x)=-\frac{1}{2\pi}\log|x|$. 
With this notation, the equivalent integral (mild) formulation to (PE), called also the Duhamel formula, reads
\begin{equation}
 \label{IPE}
u(t)=e^{t\Delta}u_0-B_0(u,u).
\end{equation}

We begin by establishing the following simple result
\begin{theorem}
 \label{theorem1}
There exist  two absolute constants $\epsilon,\beta>0$ with the following property.
Let $u_0\in E$ be such that 
$\|u_0\|_E<\epsilon$.
Then there exists a unique (mild) solution $u\in \X$ of (PE) such that
$\|u\|_\X\le \epsilon\beta$.
\end{theorem}
\medskip

The proof of Theorem~\ref{theorem1} will follow from a series of lemmata. 
\bigskip

First, we have the following estimate of the leading term in $\nabla \f$
\begin{lemma}
\label{lemma1}
Let $u\in \X$ and $\varphi$ be such that $\Delta \varphi+u=0$.
Then
\begin{equation}
 \label{nabp}
 \nabla \varphi(x,t)=\frac{c_0\,x}{|x|^2}\int_{|y|\le |x|/2} u(y,t)\d{y} + \mathcal{R}(x,t)
\end{equation}
with 
$c_0=-\frac{1}{2\pi}$ and 
the remainder $\mathcal{R}$ satisfying 
\begin{equation*}
 |\mathcal{R}(x,t)|\le C\|u\|_\X\left(t^{\frac12}+|x|\right)^{-1}.
\end{equation*}
\end{lemma}

\Proof
Indeed, let us represent the partial derivatives of $\f$, for $j=1,2$,  as 
$$\partial_j \varphi=\frac{c_0x_j}{|x|^2}*u\equiv I_1+I_2+I_3,$$
where
\begin{equation}
\begin{split} 
I_1&=\int_{|y|\le |x|/2}\frac{c_0(x_j-y_j)}{|x-y|^2}\, u(y,t)\d{y}.
\end{split}\nonumber
\end{equation}
The terms $I_2$ and $I_3$ are obtained by taking the integration domains  $\{|x-y|\le|x|/2\}$ and $\{|x-y|\ge|x|/2,\,|y|\ge|x|/2\}$, respectively, 
in the convolution integrals defining $\partial_j\f$.
It is straightforward to prove that $I_2$ and $I_3$ can be bounded by $C\|u\|_\X\left(t^{\frac12}+|x|\right)^{-1}$.
On the other hand, we can rewrite $I_1$ as
$$\frac{c_0x_j}{|x|^2}\int_{|y|\le |x|/2} u(y,t)\d{y}+ R_1(x,t).$$
An application of the Taylor formula shows that the above bound holds also for $R_1$.
 \qed

\medskip

We immediately deduce from \eqref{nabp} the following useful estimate
\begin{lemma}
 \label{lemma2}
Let $u\in \X$ and $\Delta \varphi +u=0$. Then
\begin{equation*}
\|\nabla \varphi(t)\|_{L^\infty}\le C\|u\|_\X \, t^{-\frac12}.
\end{equation*}
\end{lemma}
\medskip

The last lemma that we need is the following 

\begin{lemma}
\label{lemma3}
Let $u,v\in \X$. Then, for some constant $C_0$ independent of $u,v$ 
\begin{equation*}
 \|B_0(u,v)\|_\X\le C_0\|u\|_\X\|v\|_\X.
\end{equation*}
\end{lemma}

\Proof
We can assume, without any restriction, that $\|u\|_\X=\|v\|_\X=1$. Lemma~\ref{lemma2} implies
\begin{subequations}
\begin{equation}
\label{fa}
|u\nabla(-\Delta)^{-1} v|(x,t)\le C|x|^{-\frac32}t^{-\frac34},
\end{equation}
and also 
\begin{equation}
 \label{fb}
|u\nabla(-\Delta)^{-1} v|(x,t)\le C|x|^{-2}t^{-\frac12}.
\end{equation}
\end{subequations}

We denote the gradient of the heat semigroup kernel $g_t$ by 
\begin{equation}
 \label{Ga}
G(\cdot,t)\equiv\nabla_x \left(\frac{1}{4\pi t}e^{-\frac{1}{4t}|\,\cdot\,|^2}\right).
\end{equation}
Then we may represent $B_0$ as 
\begin{equation}
\begin{split} 
B_0(u,v)(x,t) &=\int_0^t \!\!\int G(x-y,t-s)(u\nabla(-\Delta)^{-1}v)(y,s)\d{y}\d{s}\\
&\equiv J_1+J_2,
\end{split}\nonumber
\end{equation}
where $J_1=\int_0^t\!\!\int_{|y|\le |x|/2}\dots$ and $J_2=\int_0^t\!\!\int_{|y|\ge |x|/2}\dots$\ .
Using the estimate 
\begin{equation*}
 \label{b1g}
|G(x-y,t-s)|\le C|x-y|^{-\frac52}(t-s)^{-\frac14}
\end{equation*}
and inequality~\eqref{fa}, we get the bound 
\begin{equation}
\label{p11}
|J_1(x,t)|\le C|x|^{-2}.
\end{equation}
Another possible estimate is
\begin{equation}
\label{p12bis}
|J_1(x,t)|\le C|x|^{-\frac32}t^{-\frac14},
\end{equation}
which is obtained using the bound
\begin{equation*}
 \label{bg2}
|G(x-y,t-s)|\le C|x-y|^{-2}(t-s)^{-\frac12}.
\end{equation*}

On the other hand, from the property
\begin{equation}
 \label{bg3}
\|G(\cdot,t-s)\|_{L^1}= c(t-s)^{-\frac12}
\end{equation}
and inequality~\eqref{fb}, we obtain
\begin{equation}\label{j2x}
 |J_2(x,t)|\le C|x|^{-2}.
\end{equation}
As before, we have also the bound 
\begin{equation}
\label{p12}
|J_2(x,t)|\le C|x|^{-\frac32}t^{-\frac14}.
\end{equation}
This second estimate is deduced from~\eqref{fa}.

Then, using \eqref{p11}, \eqref{j2x}, we obtain the space decay estimate
\begin{equation}
 \label{sdb}
|B_0(u,v)|(x,t)\le C|x|^{-2}
\end{equation}
and from \eqref{p12bis}, \eqref{p12} --- a provisory (not optimal) time decay estimate
\begin{equation*}
\|B_0(u,v)(t)\|_{L^{\frac43,\infty}}\le Ct^{-\frac14}.
\end{equation*}
But we may represent $B_0$ as 
\begin{equation}
B_0(u,v)(t)=e^{\frac{t}{2}\Delta}B_0(u,v)(t/2)+\int_{t/2}^t G(t-s)*(u\nabla(-\Delta)^{-1} v)(s)\d{s}.\nonumber
\end{equation}
Thus, applying the weak Young-type inequality for convolutions in Lorentz spaces $L^{\frac43,\infty}*L^{4,1}\subset L^\infty$, see \cite{LR}, 
and the equality obtained from the scaling laws in Lorentz spaces,  
\begin{equation}
 \left\|\frac{1}{4\pi t} e^{-\frac{1}{4t}|\,\cdot\,|^2}\right\|_{L^{4,1}}=c t^{-\frac34},
\end{equation}
we finally get
\begin{equation}
\label{tdd} 
\begin{split}
\|B_0(u,v)(t)\|_{L^\infty} &\le Ct^{-1}+C \int_{t/2}^t (t-s)^{-\frac12}\|u\nabla(-\Delta)^{-1} v(s)\|_{L^\infty}\d{s}\\
&\le Ct^{-1}.
\end{split}
\end{equation}
Combining inequalities~\eqref{sdb} and~\eqref{tdd} we get $B_0(u,v)\in \X$, together with its continuity with respect to $u$ and $v$.
\qed

\medskip

\Proof  Note that (using the duality ${\mathcal S}-{\mathcal S}'$)  we have $e^{t\Delta}u_0\to u_0$ in ${\mathcal S}'$ as $t\to 0$. 
The conclusion of Theorem~\ref{theorem1} follows in a standard way (cf., e.g., \cite{LR,B-SM,BCGK}) from the contraction fixed point theorem.
\qed

\bigskip

\section{The parabolic-parabolic system}
Let $\tau>0$ be a fixed parameter. We consider the system (PP) whose equivalent integral formulation reads
\begin{equation}
\label{dupp}
u(t)=e^{t\Delta}u_0-\int_0^t \nabla\cdot e^{(t-s)\Delta}\biggl[u(s)\frac{1}{\tau}\nabla\int_0^s e^{\frac{1}{\tau}(s-\sigma)\Delta}u(\sigma)\d{\sigma}\biggr]\d{s}.
\end{equation}
We introduce for all $\tau\ge0$ the bilinear form  $B_\tau$ (recall that~$G$ is defined by the expression~\eqref{Ga})
\begin{equation}\label{b-t}
 B_\tau(u,v)(x,t)\equiv \int_0^t\!\!\int G(x-y,t-s)\bigl(u\,W_\tau(v)\bigr)(y,s)\d{y}\d{s},
\end{equation}
where $W_\tau(v)$ is the linear operator acting on $v$ 
\begin{subequations}
\begin{equation}
\label{tw2a}
 W_\tau(v)(x,t)=\int_0^t \frac{1}{\tau}\biggl[G\Bigl(\frac{t-\sigma}{\tau}\Bigr)*v(\sigma)\biggr](x,\sigma)\d{\sigma}\ \ \ {\rm{ for\ }}\tau>0, 
\end{equation}
with a natural convention 
\begin{equation}
\label{tw2b}
W_0(v)(x,t)=\left(\nabla(-\Delta)^{-1}v\right)(x,t). 
\end{equation}
\end{subequations}
In this way, the system (PP) is also rewritten in a compact form (cf. \eqref{IPE}) as 
\begin{equation}
 \label{simd}
u=e^{t\Delta}u_0-B_\tau(u,u).
\end{equation}
We are going to solve~\eqref{simd} in the space $\X$ exactly as was in the parabolic-elliptic case. An additional estimate, however, is needed:

\begin{lemma}
\label{lemma4}
Let $u\in \X$ and $\tau>0$. Then there exists a constant~$C^*>0$, independent of $u$ and $\tau$, such that
\begin{equation}
\label{decW}
\|W_\tau(u)(t)\|_{L^\infty} \le {C^*}t^{-\frac12}\|u\|_\X.
\end{equation}
\end{lemma}

\Proof
As usual we can and do assume $\|u\|_\X=1$. Then, for all $1<p\le \infty$ we have 
$$ |u(x,\sigma)|\le |x|^{-\frac{2}{p}}\sigma^{-1+\frac{1}{p}}.$$
This implies, for  $1<p\le \infty$,
$$ \|u(\sigma)\|_{L^{p,\infty}}\le \sigma^{-1+\frac{1}{p}}.$$

Now, we represent 
$$ W_\tau(u)= I_1+I_2,$$
where $I_1=\int_0^{t/2}\dots$ and $I_2=\int_{t/2}^{t}\dots$\ .
Evidently, we obtain the bound 
\begin{equation*}
\begin{split}
\|I_1(t)\|_{L^\infty} 
&\le C\int_0^{t/2} \left\|\frac{1}{\tau}G\Bigl(\frac{t-\sigma}{\tau}\Bigr)\right\|_{L^{2,1}} \|u(\sigma)\|_{L^{2,\infty}} \d{\sigma}\\
&\le Ct^{-\frac12}.
\end{split}
\end{equation*}

For the integral  $I_2$, let us begin with a rough bound 
\begin{equation}
\begin{split}
\|I_2(t)\|_{L^\infty} 
&\le C\int_{t/2}^t \left\|\frac{1}{\tau}G\Bigl(\frac{t-\sigma}{\tau}\Bigr)\right\|_{L^1} \|u(\sigma)\|_{L^\infty} \d{\sigma}\label{i11}\\
&\le C\tau^{-\frac12}t^{-\frac12}.
\end{split}
\end{equation} 
This bound gives the required estimate, excepted when $\tau$ belongs to a~neighbourhood of the origin.
Thus, in the sequel, it is enough to consider the case $0<\tau<\frac{1}{2}$.
Now, we further decompose
$$I_2\equiv I_{2,1}+I_{2,2},$$
where
$I_{2,1}=\int_{t/2}^{t-\tau t}\dots$ and $I_{2,2}=\int_{t-\tau t}^t\dots$\ .
Next, we are going to improve \eqref{i11}  writing  
\begin{equation}
 \label{weq}
\begin{split}
 \|I_{2,1}(t)\|_{L^\infty} 
&\le C\int_{t/2}^{t-\tau t}  \left\|\frac{1}{\tau}G\Bigl(\frac{t-\sigma}{\tau}\Bigr)\right\|_{L^{3,1}}\|u(\sigma)\|_{L^{\frac32,\infty}}\d{\sigma}\\
&\le C\tau^{\frac16}\int_{t/2}^{t-\tau t}(t-\sigma)^{-\frac76}\sigma^{-\frac13}\d{\sigma}\\
&\le Ct^{-\frac12}.
\end{split}
\end{equation}

The last estimate is
\begin{equation}
 \label{weq2}
\begin{split}
 \|I_{2,2}(t)\|_{L^\infty} 
&\le C\int_{t-\tau t}^{t}  \left\|\frac{1}{\tau}G\Bigl(\frac{t-\sigma}{\tau}\Bigr)\right\|_{L^{1}} \|u(\sigma)\|_{L^{\infty}}\d{\sigma}\\
&\le C\tau^{-\frac12}\int_{t-\tau t}^{t}(t-\sigma)^{-\frac12}\sigma^{-1}\d{\sigma}\\
&\le Ct^{-\frac12}.  
\end{split}
\end{equation}
The conclusion of Lemma~\ref{lemma4} follows from \eqref{i11}, \eqref{weq}, \eqref{weq2}.
\qed

\medskip

Lemma~\ref{lemma4} allows us to see that, if $u,\,v\in \X$ with $\|u\|_\X=\|v\|_\X=1$, then
\begin{subequations}
\begin{equation}
\label{faa}
|u\,W_\tau(v)|(x,t)\le C|x|^{-\frac32}t^{-\frac34},
\end{equation}
and
\begin{equation}
 \label{fbb}
|u\,W_\tau(v)|(x,t)\le C|x|^{-2}t^{-\frac12},
\end{equation}
for some constant $C>0$ independent of~$\tau>0$.
\end{subequations}
These are the analogous estimates as those we obtained in the parabolic-elliptic case (see inequalities~\eqref{fa} and~\eqref{fb}).
Then, using exactly the  same arguments as in the previous section, we arrive at  the following existence result

\begin{theorem}
 \label{theorem2}
There exist  two absolute constants $\epsilon^*,\beta^*>0$ with the following property.
Let $u_0\in E$ be such that 
$\|u_0\|_E<\epsilon^*$.
Then there exists a~unique (mild) solution $u\in \X$ of (PP) such that
$\|u\|_\X\le \epsilon^*\beta^*$.
\end{theorem}
\medskip

\noindent {\bf Remark.} The case of nonzero initial data $\f(0)$ can be studied in a quite similar way. 
\bigskip

\noindent {\bf Remark.} A closer look at the proofs of estimates for $\nabla\f$ in (PE) and (PP)
reveals that the behaviour of $\nabla\f$ is a bit different in these two cases. Namely, if $0\not\equiv u$ and $u(x,t)\sim (t+|x|^2)^{-1}$
(in the sense that $c_1(t+|x|^2)^{-1}\le u(x,t)\le c_2(t+|x|^2)^{-1}$ for some $c_1,\, c_2>0$), then it follows from \eqref{nabp} that
$\nabla\f(x,t)\sim |x|^{-1}\log(t+|x|^2)$ for (PE), while $\nabla\f$ is more regular:
$|\nabla\f(x,t)|\le c\left(t^{\frac12}+|x|\right)^{-1}$ in (PP) case.
In other words, letting~$\mathcal{Y}$ be the space of functions $f=f(x,t)$ such that $f^2\in\X$,
we have $\varphi^\tau\in \mathcal{Y}$ for $\tau>0$, but $(\varphi^\tau)$ does not converge in~$\mathcal{Y}$
as $\tau\to0$.
However, such an instability does not prevent from the convergence
of the densities $u^\tau\to u$ for vanishing~$\tau$.

\bigskip 
%
%
%
%
%
%
%

\section{Study of the $\tau\searrow0$ limit}
\label{section4}

We now study the convergence as $\tau\searrow0$
of solutions $u^\tau$ of the system (PP) to the corresponding solution $u$ of (PE).
A result in this direction has been obtained recently by A.~Raczy\'nski in~\cite{AR},
who established the convergence $u^\tau\to u$ in the norm $\mathcal{Y}_\alpha$,
for $\alpha\in(1,2)$, defined as
\begin{equation}
\label{ARN}
\|u\|_{\mathcal{Y}_\alpha}= \hbox{ess\!}\sup_{\!\!\!\!\!\!\!\!\!t>0, \,\xi\in\R^2} \left(1+t^{\frac{1}{2}}|\xi|\right)^{\alpha}|\widehat u(\xi,t)|.
\end{equation}
We will obtain in subsection~\ref{subsec42} a similar result using the $\X$-norm.

\subsection{Regularity properties of solutions of (PE)}

In this subsection we prepare some preliminary material.
The first Proposition consists of a regularity result with respect to the space variable for solutions
of~(PE). The second Proposition describes their regularity properties with respect to the
time  variable.

\begin{proposition}
 \label{prop1}
For all $r\in(1,2)$ there exists a constant $\epsilon_r$, with $0<\epsilon_r\le \epsilon$
(the absolute constant of~Theorem~\ref{theorem1})
such that, if $\|u_0\|_{E}<\epsilon_r$, then the solution of (PE) constructed
in~Theorem~\ref{theorem1} satisfies
\begin{align}
\label{nabu}
 \left\|\nabla u(t)\right\|_{L^{r,\infty}}\le C t^{-\frac{3}{2}+\frac{1}{r}},
\end{align}
for some constant $C=C(u_0,r)$ independent~of $t$.
\end{proposition}

\Proof
We use a standard argument involving the subspace $\X_r\subset \X$ defined by
\begin{align*}
 \X_r=\Bigl\{u\in \X,\;
\exists\,C\colon
\|\nabla u(t)\|_{L^{r,\infty}}\le C t^{-\frac{3}{2}+\frac{1}{r}} \Bigr\},
\end{align*}
and equipped with its natural norm.
Recalling that the kernel $\frac{c_0x}{|x|^2}$ of the operator
$\nabla(-\Delta)^{-1}$ belongs to $L^{2,\infty}$,
first we deduce from the Young inequality 
\begin{align}
\label{nabbu}
\left\|\nabla^2 (-\Delta)^{-1} v(t)\right\|_{L^{\alpha,\infty}} 
\le C t^{-\frac{3}{2}+\frac{1}{r}} \|v\|_{\X_r},
\qquad
\textstyle\frac{1}{\alpha}=\frac{1}{r}-\frac{1}{2}.
\end{align}
Next, from the H\"older inequality (noticing that $\|u(t)\|_{L^{2,\infty}}\le Ct^{-\frac12}\|u\|_{\X}$),
\begin{align}
 \label{nabbbu}
\|u\nabla^2(-\Delta)^{-1}v(t)\|_{L^{r,\infty}}\le Ct^{-2+\frac1r}\|u\|_{\X}\|v\|_{\X_r}.
\end{align}
The generalization of the classical inequalities  to Lorentz spaces can be found,  e.g., in~\cite{LR}.

We claim that the bilinear operator~$B_0$ introduced in~\eqref{B0}
is boundedly defined: $B_0\colon\X_r\times\X_{r}\to  \X_r$.
Indeed, for $\|u\|_{\X_r}=\|v\|_{\X_r}=1$,
we combine the estimates 
$$\|\nabla G(t-s)\|_1\le C(t-s)^{-1}$$ 
and \eqref{bg3} 
with the inequality (a consequence of $u\in\X$) 
\be
\|u(s)\|_{L^{r,\infty}}\le Cs^{-1+\frac1r},\label{Lor}
\ee
the estimate $\|\nabla(-\Delta)^{-1}v(s)\|_\infty\le Cs^{-\frac12}$ obtained from~Lemma~\ref{lemma2},
and~\eqref{nabbbu}.
Then we arrive at 
\begin{align*}
\lefteqn{\|\nabla B_0(u,v)(t)\|_{L^{r,\infty}} }\\
&\qquad
\le C \int_0^{t/2} (t-s)^{-1}\|u(s)\|_{L^{r,\infty}}\|\nabla(-\Delta)^{-1}v(s)\|_\infty\,{\rm d}s\\
&\qquad\qquad
+C\int_{t/2}^t (t-s)^{-\frac12} \Bigl(\|\nabla u(s)\|_{L^{r,\infty}} \|\nabla(-\Delta)^{-1}v(s)\|_\infty
  +s^{-2+\frac{1}{r}}\Bigr)\,{\rm d}s\\
&\qquad
\le Ct^{-\frac{3}{2}+\frac{1}{r}}.
\end{align*}

Moreover, for $u_0\in E$, we have
$ |e^{\frac{t}{2}\Delta}u_0(x)|\le C(t+|x|^2)^{-1}$.
Hence,
$ \|e^{\frac{t}{2}\Delta}u_0\|_{L^{\gamma,\infty}}\le C t^{-1+\frac{1}{\gamma}}$ for $1<\gamma<\infty$.
We now choose  $\beta,\gamma\in(1,\infty)$ such that $1+\frac{1}{r}=\frac{1}{\beta}+\frac{1}{\gamma}$.
Then, the semigroup property of the heat kernel $g_t$, and the fact that $\nabla g_{t/2}\in L^{\beta,1}$, imply
\begin{align}
\|\nabla e^{t\Delta}u_0\|_{L^{r,\infty}} \le Ct^{-\frac{3}{2}+\frac{1}{r}}.\nonumber
\end{align}

\smallskip
\noindent Now the usual the application of the contraction mapping theorem,
in a~closed ball of small radius in the space $\X_r$, allows us to conclude.
\qed

The following proposition is the first crucial tool for our stability result.
It provides the H\"older regularity, with respect to the time variable,
of solutions of (PE) in Lorentz spaces.

\begin{proposition}
 \label{prop2}
Let $1<r<2$ and $u_0\in E$, such that $\|u_0\|_E<\epsilon_r$.
Then the solution~$u$ of (PE) constructed in Proposition~\ref{prop1} satisfies for all $0<t'<t$ 
\begin{equation}
 \label{lipt}
 \|u(t)-u(t')\|_{L^{r,\infty}}\le C(t-t')^{\frac12}(t')^{-\frac{3}{2}+\frac{1}{r}},
\end{equation}
for some $C=C(u_0,r)$ independent of $t$ and $t'$.
\end{proposition}

\Proof
It is enough to show that both $e^{t\Delta}u_0-e^{t'\Delta}u_0$ and 
$B_0(u,u)(t)-B_0(u,u)(t')$ satisfy the required bound in the $L^{r,\infty}$-norm.

From the identity 
\begin{equation*}
 \begin{split}
  e^{t\Delta}u_0(x)-e^{t'\Delta}u_0(x)
&=\int \Bigl( e^{t'\Delta}u_0(x-y)-e^{t'\Delta}u_0(x)\Bigr) g_{t-t'}(y)\,{\rm d}y\\
&=-\int\!\!\int_0^1\nabla e^{t'\Delta}u_0(x-\theta y)\cdot y g_{t-t'}(y)\,{\rm d}y\,{\rm d}\theta,
 \end{split}
\end{equation*}

we get
\begin{equation*}
\begin{split} 
 \|e^{t\Delta}u_0-e^{t'\Delta}u_0\|_{L^{r,\infty}} 
&\le 
C \|\nabla e^{t'\Delta}u_0\|_{L^{r,\infty}} \|y g_{t-t'}\|_{1}\\
&\le C(t-t')^{\frac12}(t')^{-\frac{3}{2}+\frac{1}{r}}.
\end{split}
\end{equation*}


Now, we can write
$$ B_0(u,u)(t)-B_0(u,u)(t')=A_1+A_2,$$
with
$$A_1\equiv\int_0^{t'}\bigl( G(t-s)-G(t'-s)\bigr)*\left(u\nabla(-\Delta)^{-1}u\right)(s)\,{\rm d}s$$
and
$$A_2\equiv\int_{t'}^t G(t-s)*\left(u\nabla(-\Delta)^{-1}u\right)(s)\,{\rm d}s.$$
Recall that from $u\in\X$ we deduce \eqref{Lor}. 
Combining this with the estimate of~Lemma~\ref{lemma2} we get
\begin{equation}
 \label{ineN}
 \left\|u\nabla(-\Delta)^{-1}u(s)\right\|_{L^{r,\infty}}\le Cs^{-\frac{3}{2}+\frac{1}{r}}.
\end{equation}
This immediately yields
\begin{equation*}
 \begin{split}
  \|A_2\|_{L^{r,\infty}}
\le \int_{t'}^t (t-s)^{-\frac12}s^{-\frac{3}{2}+\frac{1}{r}}\,{\rm d}s
\le C(t-t')^{\frac12}\bigl(t'\bigr)^{-\frac{3}{2}+\frac{1}{r}}.
 \end{split}
\end{equation*}

The estimate of $A_1$ is slightly more involved.
We start with the identity
$$ G(t-s)-G(t'-s)=\left(e^{(t-t')\Delta}-{\rm Id}\right)G(t'-s).$$
The action of the convolution operator with the function on the right hand side is studied
via the following variant of a result established in~\cite{Mey99}.
\begin{lemma}
 \label{lemma5}
Denote the Calder\'on operator by $\Lambda=(-\Delta)^{\frac12}$.
Then, for some constant $C$ depending only on~$r\in(1,\infty)$ 
\begin{equation*}
 \bigl\| \bigl(e^{t\Delta}-{\rm Id}\bigr)f\bigr\|_{L^{r,\infty}}\le C t^{\frac12}\|\Lambda f\|_{L^{r,\infty}}.
\end{equation*}
\end{lemma}

\Proof
Writing $f=\Lambda^{-1}\Lambda f$, we see that
$$\bigl(e^{t\Delta}-{\rm Id}\bigr) f=\Phi_t*(\Lambda f),$$
where,
$$\widehat \Phi_t(\xi)= t^{\frac12}\widehat \Phi(t^{\frac12}\xi),$$
and
$$\widehat \Phi(\xi)=(e^{-|\xi|^2}-1)|\xi|^{-1}.$$
It only remains to show that $\Phi\in L^1(\R^2)$,
which is immediate. Indeed, it is well known, and easy to check,
that
$\widehat\Psi(\xi)=|\xi|e^{-|\xi|^{2}}$ defines a function $\Psi\in L^1(\R^2)$
(for example, with the method described in~\cite{FuT}, one obtains $|\Psi(x)|\le C(1+|x|)^{-3}$
and $|\nabla\Psi(x)|\le C(1+|x|)^{-4}$).
We conclude applying the Bochner inequality to the identity
$$\Phi(x)=-2\int_1^\infty \Psi(\eta\,x)\,{\rm d}\eta.$$
\qed

Using this Lemma we deduce
\begin{equation*}
\begin{split}
 \|A_{1}\|_{L^{r,\infty}}
&\le 
C(t-t')^{\frac12} \int_0^{t'} \|\Lambda\nabla g_{t'-s}*\left(u\nabla(-\Delta)^{-1}u\right)(s)\|_{L^{r,\infty}}\,{\rm d}s\\
&= C(t-t')^{\frac12}(A_{1,1}+A_{1,2}),
\end{split}
\end{equation*}
where $A_{1,1}$ and $A_{1,2}$ are obtained splitting the integral at $s=t'/2$.

But, as the function $\Psi$ introduced in the proof of Lemma~\ref{lemma5} satisfies $\nabla\Psi\in L^1(\R^2)$,
we see by a simple rescaling that
$$\|\Lambda\nabla g_{t'-s}\|_1\le C(t'-s)^{-1}.$$
Combining this estimate with inequality~\eqref{ineN}, we get
$$ A_{1,1} \le C(t')^{-\frac{3}{2}+\frac{1}{r}},
\qquad 1<r<2.$$

For treating~$A_{1,2}$, we combine the estimate
$$\|\Lambda g_{t'-s}\|_1\le C(t'-s)^{-\frac12}$$
with the inequality (for $1<r<2$)
\begin{equation*}
\|(\nabla u)\left(\nabla(-\Delta)^{-1}u\right)(s)\|_{L^{r,\infty}}
+\|u\left(\nabla^2(-\Delta)^{-1}u\right)(s)\|_{L^{r,\infty}}
\le Cs^{-2+\frac{1}{r}},
\end{equation*} 
obtained by applying~\eqref{nabu}, Lemma~\ref{lemma2}, and \eqref{nabbbu} with $u=v$. 
We get as before 
$$ A_{1,2} \le C (t')^{-\frac{3}{2}+\frac{1}{r}},
\qquad 1<r<2,$$
and this concludes the proof of Proposition~\ref{prop2}.
\qed
%
%
%
%
%
%

\subsection{The vanishing~$\tau$ limit}
\label{subsec42}

After Proposition~\ref{prop2}, the second crucial step 
for the study of the limit as $\tau\searrow0$ consists in  the asymptotic
analysis of the linear operators $W_\tau$, $\tau\ge 0$, 
introduced in \eqref{tw2a}--\eqref{tw2b}. This is the purpose of the following lemma.

\begin{lemma} \label{lemma6}
Let $\varepsilon=\varepsilon(\tau)$ be an arbitrary function, strictly increasing and continuous on
$[0,1]$, such that $\varepsilon(0)=0$.
Let also $1<r<2$ and~$u$ be a~function satisfying, for $0<t'<t$,
\begin{equation}
 \label{hypl6}
\begin{split}
&\|u(t)\|_{L^{r,\infty}}\le Ct^{-1+\frac{1}{r}},\\
&\|u(t)-u(t')\|_{L^{r,\infty}}\le C(t-t')^{\frac12}(t')^{-\frac{3}{2}+\frac{1}{r}},  
\end{split}
\end{equation}
with a constant $C$ independent of $t,t'$.
Then, for all $t>0$, $\tau\in[0,1]$, and 
for another constant $C$, independent of~$t$ and~$\tau$,
\begin{equation}
 \label{asrw}
\varepsilon(\tau)\sup_{t>0}
t^{\frac12}\left\|\left(W_\tau(u)-W_0(u)\right)(t)\right\|_\infty
\le C\tau^{\frac{1}{r}-\frac{1}{2}}.
\end{equation}
In particular, if $u_0\in E$ is small enough (for example, $\|u_0\|_{E}\le \epsilon_{3/2}$), 
then the corresponding solution~$u$  of (PE) constructed in Proposition~\ref{prop1} satisfies 
\begin{equation}
\label{limt}
\lim_{\tau\to0} \,\,\sup_{t>0} t^{\frac12}\left\|\left(W_\tau(u)-W_0(u)\right)(t)\right\|_\infty=0.
\end{equation}
\end{lemma}

\Proof
Without any restriction we can assume that $0<\varepsilon(\tau)<\frac{1}{2}$ for positive~$\tau$.
Define $\tilde\varepsilon(\tau)$ such that $\varepsilon=\tilde\varepsilon^{\frac{1}{r}-\frac{1}{2}}$.
Borrowing from \cite{AR} the idea of splitting the time integral
using intervals depending on~$\tau$,
we write
$$ W_\tau(u)-W_0(u)\equiv J_1+J_2+J_3,$$
where
$$ J_1(t)=\int_0^{t(1-\tilde\varepsilon(\tau))}\biggl[\frac{1}{\tau}G\Bigl(\frac{t-s}{\tau}\Bigr)*u(s)\biggr]\,{\rm d}s,$$
next
$$J_2(t)= \int_{t(1-\tilde\varepsilon(\tau))}^t \frac{1}{\tau}G\Bigl(\frac{t-s}{\tau}\Bigr)*u(t)\,{\rm d}s
\; -\;W_0(u)(t),$$
and
$$J_3(t)= \int_{t(1-\tilde\varepsilon(\tau))}^t \frac{1}{\tau}G\Bigl(\frac{t-s}{\tau}\Bigr)*\bigl[u(s)-u(t)\bigr]\,{\rm d}s.$$

From the first relation of~\eqref{hypl6} and the Young inequality in Lorentz space
(using that, by~\eqref{Ga}, $G(\cdot, t)\in L^{r',1}(\R^2)$, where $r'$ is the conjugate exponent),
we get
\begin{equation*}
 \begin{split}
  \|J_1(t)\|_\infty 
&\le C\tau^{\frac{1}{r}-\frac{1}{2}}\int_0^{t(1-\tilde\varepsilon(\tau))}(t-s)^{-\frac{1}{2}-\frac{1}{r}}s^{-1+\frac1r}\,{\rm d}s\\
&\le C\left(\frac{\tau}{\tilde\varepsilon(\tau)}\right)^{\frac{1}{r}-\frac{1}{2}}t^{-\frac12}.
 \end{split}
\end{equation*}
Notice that this estimate of~$J_1$ is exactly what we need for~\eqref{asrw}.

As for $J_2$, we see from a simple computation via the Fourier transform that
$$ J_2(t)=-g_{t\tilde\varepsilon(\tau)/\tau}*W_0(u)(t).$$
If $\frac{1}{\alpha}=\frac{1}{r}-\frac{1}{2}$,
then we deduce from the usual weak-convolution estimates that  $W_0(u)=\nabla(-\Delta)^{-1}u$ is bounded
in the $L^{\alpha,\infty}$-norm, by $Ct^{-1+\frac1r}$.
Applying once more the Young inequality (using now $g_{t\tilde\varepsilon(\tau)/\tau}\in L^{\alpha',1}$),
we get, as before,
$$ \|J_2(t)\|_\infty\le C\left(\frac{\tau}{\tilde\varepsilon(\tau)}\right)^{\frac{1}{r}-\frac{1}{2}}t^{-\frac12}.$$

Applying the second of inequalities~\eqref{hypl6}, we obtain immediately
$$ \|J_3(t)\|_\infty \le C\tau^{\frac{1}{r}-\frac{1}{2}}t^{-\frac12},$$
which is even better than what we need.
This proves the inequality~\eqref{asrw}.
Choosing, for example, $r=\frac{3}{2}$ and $\varepsilon(\tau)=\frac{1}{2}\tau^{1/12}$
proves the last claim \eqref{limt} of Lemma~\ref{lemma6}.
\qed

We are now in the position of establishing our first main result 

\begin{theorem}
\label{theo3}
There exists an absolute constant $\epsilon'>0$ (a priori smaller than the constants $\epsilon,\epsilon^*>0$
in Theorems~\ref{theorem1} and~\ref{theorem2}),
such that if $u_0\in E$, $\|u_0\|_E<\epsilon'$, 
then denoting by $u\in \X$ the solution of (PE)
and  $u^\tau\in\X$ the solution of (PP) constructed in the previous theorems,
we have  as $\tau\searrow0$ 
\begin{equation*}
 u^\tau\to u \qquad\hbox{in\ \  $\X$}.
\end{equation*}
\end{theorem}

\Proof
The proof follows easily from Lemma~\ref{lemma6}.
Indeed, from the integral equations~\eqref{IPE} and~\eqref{simd},
the bilinearity of~$B_\tau$ and $B_0$, and the smallness of the solutions
$u^\tau$ and $u$, we have 
(similarly as in~\cite{AR}, where two terms in the bilinear expansion can be absorbed by the
left hand side) 
$$\|u^\tau-u\|_{\X}\le C\|B_\tau(u,u)-B_0(u,u)\|_{\X}.$$

But, by the definition of~$W_\tau$ and $W_0$ (see \eqref{tw2a}-\eqref{tw2b}),
\begin{equation*}
B_\tau(u,u)(t)-B_0(u,u)(t)=\int_0^t G(t-s)*\bigl(u(W_\tau(u)-W_0(u))\bigr)(s)\,ds.
\end{equation*}
Argueing as in the proof of Lemma~\ref{lemma3},
we obtain
\begin{equation*} 
\|B_\tau(u,u)-B_0(u,u)\|_{\X}\le C\|u\|_{\X}\Bigl(\sup_{t>0}{t^{\frac12}}\|W_\tau(u)-W_0(u)(t)\|_{\infty}\Bigr).
\end{equation*}
If $\epsilon'>0$ is small enough, then Lemma~\ref{lemma6} can be applied to
the solution~$u$ of (PE),
implying that the right hand side of the above inequality has a vanishing limit for small~$\tau$.
This finally gives
$$\|u^\tau-u\|_{\X}\to0 \qquad\hbox{as \ \ $\tau\searrow 0$}.$$
\qed

\begin{remark}
\begin{rm}
Notice that, the smaller the norm $\|u_0\|_E$, the faster the convergence $u^\tau\to u$ as \hbox{$\tau\to0$}.
This is due to the fact that for very small data it is possible to
apply Lemma~\ref{lemma6} with~$r$ close to~$1$
(despite the constants in our estimates blow up as \hbox{$r\searrow1$}).
More precisely, our arguments show that for any $0<\delta<\frac{1}{2}$,
one can find a constant $C>0$ and $\epsilon(\delta)>0$ such that,
for $\|u_0\|_E\le \epsilon(\delta)$, one has $\|u^\tau-u\|_{\X}\le C\tau^{\frac{1}{2}-\delta}$
for all \ $0\le\tau\le1$.
\end{rm}
\end{remark}

%
%
%

\subsection{The parabolic-elliptic limit in stronger topologies}
If $u_0\ge0$ is small in  the $E$-norm, and 
belongs to a smaller space, for example, $u_0\in E\cap L^1$,
then the solutions $u^\tau$ and $u$ of (PP) and (PE) will
remain in $L^1$, uniformly in time, during their evolution.
Hence, it is natural to ask whether the convergence $u^\tau\to u$  holds
also in the natural norm of $L^\infty((0,\infty);L^1)$.
Our next theorem provides a positive answer.
As the proof of this fact does not really depend
on a particular topology under consideration,
it seems appropriate to consider a more abstract setting.
\medskip

We denote by $\mathcal{L}$ any
\emph{shift invariant Banach space of local measures\/},  see \cite[Ch. 4]{LR} for their definition and main
properties.
These are Banach spaces of distributions, continuously embedded
in $\mathcal{D}'(\R^2)$.
Moreover, they are known to satisfy the following properties
(for some constant $C>0$ depending only on~$\mathcal{L}$),
\begin{enumerate}
 \item For all $f\in\mathcal{L},\;g\in L^1(\R^2)$, the convolution product $f*g$ is well defined in~$\mathcal{L}$
and $\|f*g\|_{\mathcal{L}}\le C\|f\|_{\mathcal{L}}\|g\|_1$.
\item For all $f\in\mathcal{L},\;h\in L^\infty(\R^2)$, the pointwise product $fh$ is well defined in~$\mathcal{L}$
and $\|fh\|_{\mathcal{L}}\le C\|f\|_{\mathcal{L}}\|h\|_\infty$.
\item
Each bounded sequence $\{f_k\}\subset\mathcal{L}$ has a subsequence convergent in~$\mathcal{L}$,
in the distributional sense.
\end{enumerate}

\noindent
Obvious examples of spaces satisfying these properties 
(and which are indeed shift invariant space of local measures)
are the $L^p$-spaces, $1<p\le\infty$,
the Lorentz spaces $L^{p,q}$, $1<p<\infty$, $1<q\le \infty$
and the space of bounded Borel measures~$M(\R^2)=\mathcal{C}_0(\R^2)^*$.
In the latter case, such duality relations ensures Property~3.
Other interesting examples include the Morrey--Campanato spaces~$M^{p,q}$,
($1<p\le q<\infty$) and suitable multiplier spaces, see  \cite[Ch.~17]{LR}.

On the other hand, the space of pseudomeasures, 
 i.e.,  the space of tempered distributions~$f$ such that $\widehat f\in L^\infty$)
does not fulfill the second requirement.
Therefore, the stability result in the pseudomeasure topology will not be encompassed
by our next Theorem, but requires a specific (and more involved) treatment,  see~\cite{AR}.

Notice that, because of the conservation of the total mass 
for positive solutions of (PE) and (PP), the $L^1$-norm remains constant
during the evolution.
This observation will allow us to handle the case of data $u_0\in E\cap L^1$,
despite Property~3 breaks down for~$L^1$.

\begin{theorem}
\label{theo4}
Let $u_0\in E\cap \mathcal{L}$, where $\mathcal{L}$ is either a shift invariant Banach space
of local measures, or  $\mathcal{L}=L^1$.
In the latter case we require either $u_0\ge0$ or, for signed~$u_0$, that~$\|u_0\|_1$ is sufficiently small.
Then there exists a~positive constant $\tilde\epsilon$\,,
depending only on~$\mathcal{L}$, such that, if 
$$\|u_0\|_E<\tilde \epsilon,$$
then the systems (PE) and (PP) possess  unique solutions~$u$ and $u^\tau$, respectively,
such that for an absolute constant $\tilde\beta>0$,
$\|u\|_{\X} \le \tilde\beta\tilde\epsilon$ and 
$\|u^\tau\|_{\X} \le \tilde\beta\tilde \epsilon$.
In~addition,
\begin{equation*}
\sup_{t>0}\|u(t)\|_{\mathcal{L}}<\infty
\qquad\hbox{and}\qquad
\sup_{t>0}\|u^\tau(t)\|_{\mathcal{L}}<\infty.
\end{equation*}
Moreover, the conclusion
$\|u^\tau(t)-u(t)\|_{\mathcal X}\to0$
of Theorem~\ref{theo3} is strengthened 
\begin{equation*}
\sup_{t>0}\|u^\tau(t)-u(t)\|_{\mathcal{L}}\,\to0
\qquad\hbox{as \ $\tau\searrow 0$}.
\end{equation*}
\end{theorem}

\Proof
Obviously we have, for some constant $C_0>0$ independent of~$t$ 
$$ \|e^{t\Delta}u_0\|_{\mathcal{L}}\le C_0.$$
Moreover, for each $\tau\ge 0$ (we include in this way the analysis of (PE)),
we have the estimate
\begin{equation}
\label{yel}
   \|B_\tau(u,v)(t)\|_{\mathcal{L}} 
   \le \tilde C 
   \Bigl(\sup_{t>0} \|u(t)\|_{\mathcal{L}}\Bigr) \|v\|_{\X},
\end{equation}
where $\tilde C>0$ depends only on~$\mathcal{L}$.
This follows from \eqref{b-t} written as 
$$ B_\tau(u,v)(t)=\int_0^t G(t-s)*\bigl(uW_\tau(v)\bigr)(s)\,{\rm d}s$$
for each $\tau\ge0$, with the convention \eqref{tw2b}.
Therefore, using the usual estimate \eqref{bg3}, 
we have
\begin{equation*}
 \begin{split}
 \|B_\tau(u,v)(t)\|_{\mathcal{L}}\le
  C \Bigl(\sup_{t>0}\|u(t)\|_{\mathcal{L}}\Bigr)
\Bigl(\sup_{t>0}t^{\frac12}\|W_{\tau}(v)(t)\|_\infty\Bigr).
 \end{split}
\end{equation*}
The last factor is bounded by $\|v\|_{\X}$ owing to Lemma~\ref{lemma2} 
in the case~$\tau=0$, and to Lemma~\ref{lemma3} for $\tau>0$.
This yields~\eqref{yel}.

\medskip
Now, we can consider, for $\tau\ge0$, the sequence of approximating solutions 
$$ u_k^\tau=e^{t\Delta}u_0  -  B_\tau(u_{k-1}^\tau,u_{k-1}^\tau), \qquad
 k=1,2,\dots\ .$$
When $\tilde\epsilon<\min\{\epsilon,\epsilon^*\}$,  we know
by the proofs of Theorems~\ref{theorem1} and~\ref{theorem2}
that the sequence $u^\tau_k$ converges in~$\X$ to the solution~$u^\tau$ of (PE)
or (PP). Here, of course, $u=u^0$ for the solutions of (PE).

On the other hand,  applying recursively~\eqref{yel}, we get
$u^\tau_k(t)\in\mathcal{L}$ for all $k$ and
$$ \sup_{t>0}\|u^\tau_k(y)\|_{\mathcal{L}}\le C_0
     +\tilde\beta\tilde C\tilde\epsilon
     \Bigl(\sup_{t>0}\|u^\tau_{k-1}(y)\|_{\mathcal{L}}\Bigr),$$
with $\tilde\beta=\max\{\beta,\beta^*\}$ (the constants obtained in Theorems~\ref{theorem1} and~\ref{theorem2}).
Iterating this inequality we arrive at 
$$ \sup_{t>0}\|u^\tau_k(t)\|_{\mathcal{L}}\le C'<\infty,$$
provided $\tilde\beta\tilde C\tilde\epsilon<1$.
If~$\mathcal{L}$ is a shift invariant Banach space of local measures, from Property~3
we get  for all $\tau\ge0$ 
$$ \sup_{t>0}\|u^\tau(t)\|_{\mathcal{L}}\le C'<\infty,$$
where $C'>0$ is independent on $\tau$.
Of course,
the last claim remains valid
in the case $\mathcal{L}=L^1$ and $u_0\ge0$
(notice that the smallness of~$\|u_0\|_E$ prevents
blow up results that could occur, otherwise,
when the second moment of~$u_0$ are finite
and $\int u_0>8\pi$. See the introduction and the references therein quoted).
If we remove the assumption~$u_0\ge0$, we can obtain the same conclusion
provided $\|u_0\|_1$ is sufficiently small.
Indeed, we see from inequality~\eqref{yel}
that the fixed point argument applies in the space~$L^\infty_t(L^1)\cap\mathcal{X}$.
\medskip

We now discuss the stability,
including also the case $\mathcal{L}=L^1$.
From the bilinearity of $B_\tau$, the mixed estimate~\eqref{yel}
and the smallness of the solutions $u$ and $u^\tau$
(allowing two terms of the bilinear expansions 
to be absorbed by the left hand side), we obtain 
\begin{equation*}
 \begin{split}
  \|u^\tau(t)-u(t)\|_{\mathcal{L}}
&\le C\sup_{t>0} \|B_\tau(u,u)-B_0(u,u)(t)\|_{\mathcal{L}}.
 \end{split}
\end{equation*}
Argueing as at the end of the proof of Theorem~\ref{theo3}, we arrive at 
$$ \|B_\tau(u,u)-B_0(u,u)(t)\|_{\mathcal{L}}\le C\Bigl(\sup_{t>0}\|u(t)\|_{\mathcal{L}}\Bigr)
    \Bigl(\sup_{t>0}  t^{\frac12}\|(W_\tau(u)-W_0(u))(t)\|_\infty\Bigr),$$
and the conclusion follows from Lemma~\ref{lemma6}.
\qed
\medskip

\noindent {\bf Remark.} 
As an application of this general result, let us observe that 
taking $\mathcal{L}=L^\infty$, we obtain for $u_0\in E\cap L^\infty$, with $u_0$ small in the 
$E$-norm,
$$  u^{\tau}\to u$$
as $\tau\searrow 0$, uniformly in $(x,t)\in\R^2\times[0,\infty)$.

%
%
%
%
%

\section{Blow up of complex valued solutions of the parabolic-parabolic system}

Consider the system (PP) in the   space $\R^d$ with any $d\ge 1$.   
Passing to the Fourier transforms, we may write the Duhamel formula \eqref{dupp}  in the form 
\begin{equation}\label{F}
\begin{split}
&\widehat u(\xi,t)=e^{-t|\xi|^2}\widehat u_0(\xi)\\
 &\qquad\qquad+\int_0^t\!\!\int_0^s\!\!\int_{\R^d} 
\frac{\xi\cdot\eta}{\tau} 
e^{-(t-s)|\xi|^2}e^{-\frac{1}{\tau}(s-\sigma)|\eta|^2}
\widehat u(\xi-\eta,s)\widehat u(\eta,\sigma)\d{\eta}\d{\sigma}\d{s}.
\end{split}
\end{equation}

Our goal is to construct a class of complex valued initial data, such that the corresponding solutions blow
up in finite time, in any classical norm.
For~$a\in\R$, we denote by $\dot B^{a,\infty}_\infty$ the homogeneous Besov space, which can also
be identified with the H\"older--Zygmund space~$\dot{\mathcal{C}}^a$.
As it is well known (see~\cite{LR, MS}), most of the classical functional spaces (including all homogeneous
Triebel--Lizorkin, and thus Lebesgue or Sobolev spaces)
are continuously embedded in~$\dot B^{a,\infty}_\infty$ for some real~$a$.

\begin{theorem} 
\label{theoremCC}
There exists  $w_0\in\S(\R^d)$ such that $\widehat{w}_0(\xi)\ge 0$, $\|\widehat{w}_0\|_{L^1}=1$.
 If $A$ is sufficiently large, then any (tempered) distributional solution of \eqref{F} with $u_0=Aw_0$ \
(and thus any solution of (PP)) blows up in a finite time.
More precisely, for some time $t^*<\infty$ and each $a\in\R$,\  $\|u(t^*)\|_{\dot B^{a,\infty}_\infty}=\infty$
holds.
\end{theorem} 

\medskip

Our approach is closely related to that in \cite[Theorem 3.1]{BB} which followed the argument in \cite{MS} for the ``cheap'' Navier--Stokes equations.
We produce some estimates from below of the Fourier transform of any  solution with $u(0)=u_0$ that can be obtained via the iteration procedure for \eqref{F} with 
$$u_{k+1}=e^{t\Delta}u_0-\int_0^t \nabla\cdot e^{(t-s)\Delta}\frac{1}{\tau}\biggl[u_k(s)\nabla\int_0^s e^{\frac{1}{\tau}(s-\sigma)\Delta}u_k(\sigma)\d{\sigma}\biggr]\d{s},$$ 
for $k=1,2,\dots\ $. 
This recurrence relation, in general, does not preserve the positivity of the Fourier transform of $u_0$. 
\noindent 
This leads us to restrict our attention to data 
of the form $\widehat u_0=A\widehat w_0$ with $w_0\in\mathcal{S}(\R^d)$,
such that
 $\hbox{supp }{\widehat w}_0\subset \left\{\xi\in\R^d\colon 2^{-1}\le \xi_1\le |\xi|\le 1\right\}$.
Define, for $k=0,1,\dots$ the set
$$E_k=\left\{\xi\in\R^d\colon 2^{k-1}\le \xi_1\le |\xi|\le 2^k\right\}.$$
Then we see that for $w_k=w_0^{2^k}$, 
 we have  $\widehat w_k= (2\pi)^{-d}\widehat w_{k-1}*\widehat w_{k-1}$, and therefore, $\hbox{supp\,}\widehat w_k\subset E_k$.
This implies that if, in addition, $\widehat w_0(\xi)\ge0$, then the positivity of the Fourier transform
will be preserved by the sequence~$u_k$, and so by~$u$.
Next Lemma tells us more:

\begin{lemma}
\label{indu}
For all $k=0,1,2,\ldots$\, , we have
\begin{equation}
\label{fort}
\widehat u(\xi,t)\ge \beta_ke^{-2^k t}{\bf 1}_{t_k\le t<t^*}(t)\widehat w_k(\xi),
\end{equation}
where $\{\beta_k\}$ and $\{t_k\}$ are two sequences defined below in~\eqref{tk} and \eqref{betak}.
\end{lemma}

\Proof
For $k=0$, the conclusion immediately follows from
$$\widehat u(\xi,t)\ge Ae^{-t|\xi|^2}\widehat w_0(\xi),$$
provided we choose $\beta_0=A$ and $t_0=0$.
Let $k\ge1$. Assume that the inequality of the lemma holds for $k-1$. Then, for all $t_k\le t<t^*$,
\begin{equation*}
 \begin{split}
 \widehat u(\xi,t)
&\ge \int\limits_{t_{k-1}}^{t}\int\limits_{t_{k-1}}^s\int\limits_{\R^d} 
\frac{\xi_1\eta_1}{\tau}
e^{-(t-s)|\xi|^2}e^{-\frac{1}{\tau}(s-\sigma)|\eta|^2}
\beta_{k-1}^2 e^{-2^{k-1}s}e^{-2^{k-1}\sigma} \times\\
&\qquad\qquad\qquad\qquad\qquad\times
\frac{\widehat w_{k-1}(\xi-\eta)\widehat w_{k-1}(\eta)}{(2\pi)^d}\d{\eta}\d{\sigma}\d{s}\\
&
\ge \int\limits_{t_{k-1}}^{t}\int\limits_{\R^d} 
(s-t_{k-1})
\frac{\xi_1 2^{k-2}}{\tau}
e^{-(t-s)|\xi|^2}e^{-\frac{1}{\tau}(t^*-t_{k-1})2^{2k-2}}
\beta_{k-1}^2 e^{-2^k s}\times\\
&\qquad\qquad\qquad\qquad\qquad\times
\frac{\widehat w_{k-1}(\xi-\eta)\widehat w_{k-1}(\eta)}{(2\pi)^d}\d{\eta}\d{s}\\
\end{split}
\end{equation*}
Thus, we can bound $\widehat u(\xi,t)$ from below as follows
\begin{equation*}
\begin{split}
\widehat u(\xi&,t) \ge \int\limits_{t_{k-1}}^{t} 
(s-t_{k-1})
\frac{2^{2k-3}}{\tau}
e^{-(t-s)2^{2k}}e^{-\frac{1}{\tau}(t^*-t_{k-1})2^{2k-2}}
\beta_{k-1}^2 e^{-2^{k}s}\widehat w_{k}(\xi)\d{s}\\
&
\ge 
\biggl(\int\limits_{t_{k-1}}^{t} 
(s-t_{k-1})
e^{-(t-s)2^{2k}}
\d{s}
\biggr)
\frac{2^{2k-3}}{\tau}
e^{-\frac{1}{\tau}(t^*-t_{k-1})2^{2k-2}}
\beta_{k-1}^2 e^{-2^k t}\widehat w_{k}(\xi)\\
&
=
\biggl(
\frac{t-t_{k-1}}{2^{2k}}-\frac{1-e^{-(t-t_{k-1})2^{2k}}}{2^{4k}}
\biggr)
\frac{2^{2k-3}}{\tau}
e^{-\frac{1}{\tau}(t^*-t_{k-1})2^{2k-2}}
\beta_{k-1}^2 e^{-2^{k}t}\widehat w_{k}(\xi)\\
&\ge
\biggl(
(t_k-t_{k-1})-\frac{1-e^{-(t^*-t_{k-1})2^{2k}}}{2^{2k}}
\biggr)
\frac{2^{-3}}{\tau}
e^{-\frac{1}{\tau}(t^*-t_{k-1})2^{2k-2}}
\beta_{k-1}^2 e^{-2^{k}t}\widehat w_{k}(\xi).
\end{split}
\end{equation*}
This suggests us to set, for some $\delta>0$,  $t^*-t_{k-1}=\delta\tau 2^{-2k+2}$.  With this choice, putting $t_0=0$, we have
\begin{equation}
 \label{tk}
t^*=\delta\tau, \qquad t_{k}=\delta\tau(1- 2^{-2k}).
\end{equation}
Then, $t_{k}-t_{k-1}= 3\delta\tau 2^{-2k}$.
We get, for $t_k\le t<t^*$,
\begin{equation*}
 \begin{split}
  \widehat u(\xi,t)
   &\ge
    (3\delta \tau-1+ e^{-4\delta\tau}) 2^{-2k-3}\tau^{-1} e^{-\delta} \beta_{k-1}^2 e^{-2^k\,t} \widehat w_k(\xi).\\
 \end{split}
\end{equation*}
This inequality is interesting only if the right hand side is positive.
Therefore we will assume that $3\delta\tau\ge1$.
We choose $\{\beta_k\}$ in  such a way that
\begin{equation*}
\beta_0=A, \qquad \beta_{k}=(3\delta \tau-1+ e^{-4\delta\tau})
 2^{-2k-3}\tau^{-1} e^{-\delta} \beta_{k-1}^2, \qquad k=1,2,\ldots\ .
\end{equation*}
This choice leads to inequality~\eqref{fort}.

In order to compute $\beta_k$, we introduce
$M_{\delta,\tau}$, such that
$$
 2^{M_{\delta,\tau}}=(3\delta \tau-1+ e^{-4\delta\tau}) e^{-\delta}\,2^{-3}\tau^{-1}.
$$ 
Notice that we have $\beta_k=2^{M_{\delta,\tau}-2k}\beta_{k-1}^2$.
We claim that, for some $a,b,c\in\R$,
$$ \beta_k= A^{2^k}2^{a+bk+c2^k}.$$
Indeed, from an easy calculation we find $b=2$,  $a=4-M_{\delta,\tau}$ and finally $c=M_{\delta,\tau}-4$, which is needed to ensure $\beta_0=A$.
 
 Hence, we have 
\begin{equation}
\label{betak}
 \beta_k = \Bigl(A\,2^{M_{\delta,\tau}-4}\Bigr)^{2^k}\,2^{4-M_{\delta,\tau}+2k}, \qquad k=0,1,\ldots\ .
\end{equation}
\qed

We conclude that when 
$$A\ge 2^{4-M_{\delta,\tau}},$$
 we have by \eqref{betak} $\beta_k\to\infty$ and, in particular, 
$\|u(t_k)\|_{L^1}=\|\widehat u(t_k)\|_\infty\to\infty$ for $k\to\infty$.
The above size condition on~$A$ can be rewritten in an equivalent form as
 $$(3\delta \tau-1+ e^{-4\delta\tau})  A \ge 2^7e^\delta\,\tau.$$
A further analysis of the lower  bounds obtained for the Fourier transform of a~candidate solution permits us to conclude, as was in \cite{BB}, that $\|u(t^*)\|_{\dot B^{a,\infty}_\infty}=\infty$ for each $a\in \R$, so that all Besov (and also $L^p$ or Triebel--Lizorkin)
norms of $u$ blow up not later than $t^*$.
Notice that, for a blow-up at $t^*=1$, we need 
$A\ge Ce^{1/\tau}\tau$, cf. \eqref{tk}. 

\bigskip

\paragraph{Acknowledgments} This research has been supported by the European Commission Marie Curie Host Fellowship for the Transfer of Knowledge ``Harmonic Analysis, Nonlinear Analysis and Probability''  MTKD-CT-2004-013389, and by the Polish Ministry of Science (MNSzW) grant --- project N201 022 32/0902.

\bibliographystyle{amsplain}

\end{document}